\documentstyle[12pt]{article}

\makeatletter

\def\LRA#1#2{\@tempdimb=\c@enumiv\@tempdima%
   \vcenter{\offinterlineskip\halign{##\cr%
   \hfil${\scriptstyle{#1}}$\hfil\crcr%
   \hbox to \@tempdimb{\rightarrowfill}\cr%
   \noalign{\kern-1ex}%
   \hbox to \@tempdimb{\leftarrowfill}\cr%
   \hfil${\scriptstyle{#2}}$\hfil\crcr}}}
\def\RA#1{\@tempdimb=\c@enumiv\@tempdima\vbox{\offinterlineskip%
   \halign{##\cr\hfil${\scriptstyle {#1}}$\hfil\crcr%
   \hbox to \@tempdimb{\rightarrowfill}\cr}}}
\def\LA#1{\@tempdimb=\c@enumiv\@tempdima\vbox{\offinterlineskip%
   \halign{##\cr\hfil${\scriptstyle {#1}}$\hfil\crcr%
   \hbox to \@tempdimb{\leftarrowfill}\cr}}}

\def\diag{\leavevmode\bgroup\setcounter{enumiv}{1}%
   \unitlength1em \@tempdima3em \def\\{\crcr&}\vbox\bgroup%
   \def\multicolumn##1##2{\multispan##1\setcounter{enumiv}{##1}%
   \hfil{##2}\hfil\setcounter{enumiv}{1}}
   \offinterlineskip\halign\bgroup\vrule height.8em depth.7em %
   width0pt##&&\hfil${\displaystyle{##}}$\hfil\cr&}
\def\enddiag{\crcr\egroup\egroup\egroup}
\makeatother
\if@twoside \oddsidemargin -20pt \evensidemargin 20pt \marginparwidth 85pt
\else \oddsidemargin 0pt \evensidemargin 0pt
\fi
 \font\symbolfont=msbm10

\def\ddd{\hbox{\symbolfont D}}
\def\eee{\hbox{\symbolfont E}}

\def\nnn{\hbox{\symbolfont N}}

\def\ttt{\hbox{\symbolfont T}}

\def\cc{\mathop{\cal C}\nolimits}
\def\dd{\mathop{\cal D}\nolimits}
\def\ee{\mathop{\cal E}\nolimits}

\def\ii{\mathop{\cal I}\nolimits}
\def\pp{\mathop{\cal P}\nolimits}

\def\ff{\mathop{\cal F}\nolimits}

\def\ww{\mathop{\cal W}\nolimits}
\def\xx{\mathop{\cal X}\nolimits}
\def\yy{\mathop{\cal Y}\nolimits}

\def\pd{\mathop{\rm pd}\nolimits}

\def\id{\mathop{\rm id}\nolimits}

\def\spep#1{\mathop{{}^{\bullet}\strut\kern-.1em{#1}}\nolimits}

\def\ext{\mathop{\rm Ext}\nolimits}
\def\hom{\mathop{\rm Hom}\nolimits}

\def\add{\mathop{\rm add}\nolimits}

\def\tr{\mathop{\rm Tr}\nolimits}

\def\Cok{\mathop{\rm Cok}\nolimits}
\def\Ker{\mathop{\rm Ker}\nolimits}

\def\Im{\mathop{\rm Im}\nolimits}

\def\mod{\mathop{\rm mod}\nolimits}

\def\add{\mathop{\rm add}\nolimits}

\def\cc{\mathop{\cal C}\nolimits}

\def\Cok{\mathop{\rm Cok}\nolimits}
\def\Ker{\mathop{\rm Ker}\nolimits}

\def\pd{\mathop{\rm pd}\nolimits}
\def\fd{\mathop{\rm fd}\nolimits}
\def\id{\mathop{\rm id}\nolimits}

\def\fin{\mathop{\rm fin.dim}\nolimits}

\def\ext{\mathop{\rm Ext}\nolimits}

\def\tr{\mathop{\rm Tr}\nolimits}

\def\grade#1{\mathop{{\rm grade}\ #1}}
\def\sgrade#1{\mathop{{\rm s.grade}\ #1}}
\def\rgrade#1{\mathop{{\rm r.grade}\ #1}}
\def\extf#1#2{\mathop{\eee_{#1}#2}}
\def\trf#1#2{\mathop{\ttt_{#1}#2}}

\def\XA{1}
\def\XAA{1.1}
\def\XAB{1.2}
\def\XAC{1.3}
\def\XB{2}
\def\XBA{2.1}
\def\XBB{2.2}
\def\XBC{2.3}
\def\XBCA{2.3.1}
\def\XBCB{2.3.2}
\def\XBCC{2.3.3}
\def\XBD{2.4}
\def\XC{3}
\def\XCA{3.1}
\def\XCAA{3.1.1}
\def\XCB{3.2}
\def\XCBA{3.2.1}
\def\XCC{3.3}
\def\XCD{3.4}
\def\XCDA{3.4.1}
\def\XCDB{3.4.2}
\def\XCDC{3.4.3}
\def\XCE{3.5}
\def\XD{4}
\def\XDA{4.1}
\def\XDAA{4.1.1}
\def\XDB{4.2}
\def\XDC{4.3}
\def\XDD{4.4}
\def\XDE{4.5}
\def\XDF{4.6}
\def\XDG{4.7}
\def\XDH{4.8}
\def\XDI{4.9}
\def\XDJ{4.10}
\def\XE{5}
\def\XEA{5.1}
\def\XEB{5.2}
\def\XEC{5.3}
\def\XED{5.4}
\def\XEE{5.5}
\def\XEF{5.6}
\def\XEG{5.7}
\def\XEGA{5.7.1}
\def\XEGB{5.7.2}
\def\XEH{5.8}
\begin{document}
\begin{center}
{\Large \bf Auslander-type conditions and cotorsion pairs}

\vskip1cm Zhaoyong Huang{\footnote{The first author was partially
supported by Specialized Research Fund for the Doctoral Program of
Higher Education (Grant No. 20060284002) and NSF of Jiangsu Province
of China (Grant No. BK2007517).}} and Osamu Iyama$^2$

\vskip0.5cm $^1$ {\footnotesize Department of Mathematics, Nanjing
University, Nanjing 210093, People's Republic of China \\(E-mail:
huangzy@nju.edu.cn)}

$^2$ {\footnotesize Graduate School of Mathematics, Nagoya
University Chikusa-Ku, Nagoya 464-8602, Japan \\(E-mail:
iyama@math.nagoya-u.ac.jp)}
\end{center}

\begin{abstract}
We study the properties of rings satisfying Auslander-type
conditions. If an artin algebra $\Lambda$ satisfies the Auslander
condition (that is, $\Lambda$ is an $\infty$-Gorenstein artin
algebra), then we construct two kinds of subcategories which form
functorially finite cotorsion pairs.
\end{abstract}

\vspace{.5em}
Noetherian rings satisfying `Auslander-type conditions' on self-injective
resolutions can be regarded as certain non-commutative analogs of
commutative Gorenstein rings.
Such conditions, especially dominant dimension and the $n$-Gorenstein
condition, play a crucial role in representation theory and non-commutative
algebraic geometry (e.g. [A], [AR2,3], [B], [C], [EHIS], [FGR], [FI], [HN], [IS], [Iy4,6], [M], [R], [Sm], [T], [W]). They are also
interesting from the viewpoint of some unsolved homological conjectures,
e.g. the finitistic dimension conjecture, Nakayama conjecture,
Gorenstein symmetry conjecture, and so on.
It is therefore important to understand non-commutative `regular' or
`Gorenstein' rings though it is still far from realized even for the case of
finite dimensional algebras.
Recently, several authors (e.g. [Hu1,2,3], [Iy1,3]) have studied some Auslander-type
conditions, e.g. the quasi $n$-Gorenstein condition, the $(l,n)$-condition,
and so on.
This paper is devoted to enlarge our knowledge of the homological
behavior of non-commutative rings.
Especially we introduce Auslander-type conditions $G_n(k)$ and $g_n(k)$ and study their properties.

Throughout this paper, let $\Lambda$ be a left and right noetherian ring (unless
stated otherwise). We denote by
\[0 \to \Lambda \to I_0(\Lambda) \to
I_1(\Lambda) \to \cdots \to I_i(\Lambda) \to \cdots\]
the minimal injective resolution of the left $\Lambda$-module $\Lambda$. We call $\Lambda$
$n$-{\it Gorenstein} if $\fd I_i(\Lambda) \leq i$ for any $0 \leq i
\leq n-1$\footnote{Notice that there are other meanings to
  the notion of $n$-Gorenstein rings (e.g. [Iw1,2][EJ][EX]).}, and call $\Lambda$ $\infty$-{\it Gorenstein} if it is
$n$-Gorenstein for all $n$. In the latter case, $\Lambda$ is also
said to satisfy the {\it Auslander condition}. It was proved by
Auslander that the notion of $n$-Gorenstein rings is left-right
symmetric [FGR; 3.7] (see \XDB\ below).
Our aim in this paper is to generalize the Auslander condition. 
In Section 3, we introduce Auslander-type conditions $G_n(k)$ and $g_n(k)$ and study their properties. 
In section 4, we concentrate on the conditions $G_n(k)$ and $g_n(k)$ for the case $k=0,1$. We give a quick proof of well-known results on these conditions, then prove our main result.
In Section 5, we apply our results to the finitistic dimension and the $(l,n)$-condition.

\vskip1em
We denote by $\mod\Lambda$ the category of finitely
generated left $\Lambda$-modules, and by $\underline{\mod}\Lambda$
the stable category of $\Lambda$ [AB]. Put $(\ )^*:=\hom_\Lambda(\
,\Lambda)$ and
\begin{eqnarray*}
\extf{n}{}:=\ext^n_\Lambda(\ ,\Lambda):&\mod\Lambda\to\mod\Lambda^{op}&\mbox{ for }n\ge0,\\
\trf{n}{}:=\tr\circ\Omega^{n-1}:&\underline{\mod}\Lambda\to\underline{\mod}\Lambda^{op}&\mbox{ for }n>0,
\end{eqnarray*}
where $\Omega
:\underline{\mod}\Lambda\to\underline{\mod}\Lambda$ is the syzygy
functor and $\tr
:\underline{\mod}\Lambda\to\underline{\mod}\Lambda^{op}$ is the
transpose functor [AB]. Let
\begin{eqnarray*}
\grade{X}:=&\inf\{ i\ge0\ |\ \extf{i}{X}\neq0\}&\mbox{ the {\it grade}},\\
\sgrade{X}:=&\inf\{ \grade{Y}\ |\ Y\subseteq X\}&\mbox{ the {\it strong grade}},\\
\rgrade{X}:=&\inf\{ i>0\ |\ \extf{i}{X}\neq0\}&\mbox{ the {\it reduced grade}}.
\end{eqnarray*}
When $\Lambda$ is an artin algebra over $R$, we denote by $\ddd:\mod\Lambda\to\mod\Lambda^{op}$ the duality induced by the Matlis duality of $R$.
Given two homomorphisms of modules, say $f:A\to B$ and $g:B\to C$, the composition homomorphism of $f$ and $g$ is denoted by $fg:A\to C$.

\vskip.5em{\bf\XA\ Main result }

\vskip.5em{\bf\XAA\ }
For subcategories $\cc_i$ ($i=1,2$) of $\mod\Lambda$, we denote by
$\ee(\cc_1,\cc_2)$ the subcategory of $\mod\Lambda$ consisting of
$C\in\mod\Lambda$ such that there exists an exact sequence $0\to
C_2\to C\to C_1\to0$ with $C_i\in\cc_i$ ($i=1,2$). For a subcategory
$\cc$ of $\mod\Lambda$, we use $\add\cc$ to denote
the subcategory of $\mod\Lambda$ consisting of all $\Lambda$-modules
isomorphic to direct summands of finite direct sums of modules in
$\cc$.

For a left $\Lambda$-module $C$, we use $\pd
C$, $\fd C$ and $\id C$ to denote the projective dimension, flat
dimension and injective dimension of $C$, respectively. A module $C$
in $\mod\Lambda$ is called $n$-{\it torsionfree} if
$\extf{i}{\tr C}=0$ for any $1 \leq i \leq n$ [AB].
For $n,m\ge0$, we define several full subcategories of $\mod\Lambda$
as follows:
\begin{eqnarray*}
\ww_n&:=&\{C\in\mod\Lambda\ |\ \rgrade{C}>n\},\\
\ff_n&:=&\{ C\in\mod\Lambda\ |\ C\mbox{ is $n$-torsionfree}\},\\
\pp_n&:=&\{ C\in\mod\Lambda\ |\ \pd_\Lambda C<n\},\\
\ii_n&:=&\{ C\in\mod\Lambda\ |\ \id_\Lambda C<n\},\\
\xx_{n,m}&:=&\ww_n\cap\ff_m,\\
\yy_{n,m}&:=&\add\ee(\ii_m,\pp_n).
\end{eqnarray*}
For example, we have $\ww_0=\ff_0=\mod\Lambda$ and $\pp_0=\ii_0=0$.
Notice that $\ee(\ii_m,\pp_n)$ is not necessarily closed under direct summands
(e.g. take $\Lambda$ to be the path algebra of the quiver $\bullet\to\bullet
\to\bullet$ and $n=m=1$).
We denote by $\ww_n^{op}$,
$\ff_n^{op}$, $\pp_n^{op}$, $\ii_n^{op}$, $\xx_{n,m}^{op}$ and
$\yy_{n,m}^{op}$ the corresponding categories for $\Lambda^{op}$.

\vskip.5em{\bf\XAB\ }
Following Salce [Sa], we call a pair $(\xx,\yy)$ of subcategories of
$\mod\Lambda$ a {\it cotorsion pair} if
\[\xx=\{C\in\mod\Lambda\ |\ \ext^1_\Lambda(C,\yy)=0\} \mbox{ and } 
\yy=\{C\in\mod\Lambda\ |\ \ext^1_\Lambda(\xx,C)=0\}.\]
Notice that we do not assume the
vanishing of higher Ext-groups $\ext^i_\Lambda$ ($i>1$).

%
%
%
%
%
%
%
%

\vskip.5em{\bf\XAC\ }
Now we can state a main result in this paper, which we shall prove in \XDI.

\vskip1em{\bf Theorem }{\it Let $\Lambda$ be an $\infty$-Gorenstein
artin algebra. Then $(\xx_{i,j-1},\yy_{i,j})$ ($i\ge0,\ j\ge1$) and
$(\yy_{i,j},\ddd\xx_{j,i-1}^{op})$ ($i\ge1,\ j\ge0$) form cotorsion pairs.}

\vskip.5em
Moreover, we will show that they are functorially finite in the sense of \XBC\ below. For example, we have cotorsion pairs $(\ww_{i},\yy_{i,1})$ ($j:=1$)
and $(\ff_{j-1},\ii_{j})$ ($i:=0$).

\vskip.5em{\bf\XB\ Preliminaries }

\vskip.5em{\bf\XBA\ }
Let us start with the following simple observation, which
will be used frequently in this paper.

\vskip.5em{\bf Lemma }{\it Let $\yy$ be a full subcategory of $\mod\Lambda$
which is closed under extensions, and $0\to C_1\to Y\to C_0\to0$ an exact
sequence with $Y\in\yy$. If there exists an exact sequence
$0\to Y_0\to X\to C_0\to0$ with $Y_0\in\yy$, then there exist
exact sequences $0\to C_1\to Y_1\to X\to0$ with $Y_1\in\yy$.}

\vskip.5em{\sc Proof }
The middle row in the
following pull-back diagram gives our desired exact sequence:
\[\begin{diag}
&&&&0&&0\\
&&&&\uparrow&&\uparrow\\
0&\RA{}&C_1&\RA{}&Y&\RA{}&C_0&\RA{}&0\\
&&\parallel&&\uparrow&&\uparrow\\
0&\RA{}&C_1&\RA{}&Y_1&\RA{}&X&\RA{}&0\\
&&&&\uparrow&&\uparrow\\
&&&&Y_0&=&Y_0&\\
&&&&\uparrow&&\uparrow\\
&&&&0&&0\\
\end{diag}\]\vskip-3em\hfill\rule{5pt}{10pt}

\vskip.5em {\bf\XBB\ }Assume that $\cc\supset \dd$ are
full subcategories of $\mod\Lambda$ and $C\in \cc$,
$D\in\dd$. A morphism $f: D\to C$ is said to be a {\it
right} $\dd$-{\it approximation} of $C$ if $\hom_\Lambda(X,f):
\hom_{\Lambda}(X, D)
\stackrel{}{\to}\hom_{\Lambda}(X, C)\to 0$ is exact for
any $X\in\dd$. A right $\dd$-approximation $f:
D\to C$ is called {\it minimal} if an endomorphism $g: D\to D$ is an
automorphisn whenever $f=gf$. The subcategory $\dd$ is said
to be {\it contravariantly finite} in $\cc$ if any
$C\in\cc$ has a right $\dd$-approximation.
Dually, we define the notions of {\it (minimal) left}
$\dd$-{\it approximations} and {\it covariantly finite
subcategories}. The subcategory of $\cc$ is said to be {\it
functorially finite} in $\cc$ if it is both contravariantly
finite and covariantly finite in $\cc$ [AR1].

\vskip.5em{\bf\XBC\ }
Let $\Lambda$ be an artin algebra. We call a
cotorsion pair $(\xx,\yy)$ {\it functorially finite} if the following
equivalent conditions are satisfied:

(1) $\xx$ is a contravariantly finite subcategory of $\mod\Lambda$.

(2) $\yy$ is a covariantly finite subcategory of $\mod\Lambda$.

(3) For any $C\in\mod\Lambda$, there exists an exact sequence
$0\to Y\to X\stackrel{f}{\to}C\to0$ with $X\in\xx$ and $Y\in\yy$.

(4) For any $C\in\mod\Lambda$, there exists an exact sequence $0\to
C\stackrel{g}{\to}Y'\to X'\to0$ with $X'\in\xx$ and $Y'\in\yy$.

\vskip.5em{\bf\XBCA\ (Wakamatsu's Lemma) } Let $\xx$ be a subcategory
of $\mod\Lambda$ which is closed under extensions. If $0\to A\to
B\stackrel{f}{\to}C$ is exact and $f$ is a minimal right
$\xx$-approximation of $C$, then $\ext^1_\Lambda(\xx,A)=0$ holds.

\vskip.5em{\bf\XBCB\ Proof of \XBC\ }(3)$\Rightarrow$(1) $f$ is a
right $\xx$-approximation of $C$ by $\ext^1_\Lambda(\xx,Y)=0$.

(1)$\Rightarrow$(3) Let $f:X\to C$ be a minimal right
$\xx$-approximaiton of $C$. Since $\Lambda\in\xx$, $f$ is
surjective. Since $\xx$ is closed under extensions, $\Ker f\in\yy$
holds by \XBCA.

(3)$\Rightarrow$(4) (cf. [Sa]) Let $0\to C\to I\to\Omega^{-1}C\to0$ and $0\to
Y\to X\to\Omega^{-1}C\to0$ be exact sequences with $I$ the injective
envelope of $C$, $X\in\xx$ and $Y\in\yy$. Applying \XBA, we have the desired
sequence.

(2)$\Leftrightarrow$(4) and (4)$\Rightarrow$(3) can be shown dually.
\hfill\rule{5pt}{10pt}

\vskip.5em{\bf\XBCC\ }Let $\Lambda$ be an artin algebra and
$(\xx,\yy)$ a pair of subcategories of $\mod\Lambda$ which are
closed under direct summands. If $\ext^1_\Lambda(\xx,\yy)=0$ and
the conditions \XBC(3)(4) are satisfied, then $(\xx,\yy)$ is
a functorially finite cotorsion pair.

\vskip.5em{\sc Proof }By the condition \XBC(3),
$\ext^1_\Lambda(C,\yy)=0$ implies $C\in\xx$. By the condition
\XBC(4), $\ext^1_\Lambda(\xx,C)=0$ implies
$C\in\yy$.\hfill\rule{5pt}{10pt}

\vskip.5em
For a full subcategory $\cc$ of $\mod\Lambda$, we denote by
$\underline{\cc}$ the corresponding subcategory of $\underline{\mod}\Lambda$.
We denote by $\Omega^n(\mod\Lambda)$ the full subcategory of $\mod\Lambda$
consisting of $C\in\mod\Lambda$ such that there exists an exact sequence
$0\to C\to P_0\to\cdots\to P_{n-1}$ with projective $P_i$.

\vskip.5em{\bf\XBD\ }
We collect some useful results.

(1) We have the following diagram whose rows are equivalences
and columns are dualities [Iy5; 1.1.1]:
\[\begin{diag}
\underline{\ww}_n&=&\underline{\xx}_{n,0}&\RA{\Omega}&\underline{\xx}_{n,1}&\RA{\Omega}&\cdots&\RA{\Omega}&\underline{\xx}_{1,n-1}&\RA{\Omega}&\underline{\xx}_{0,n}&=&\underline{\ff}_n\\
&&\downarrow^{\tr}&&\downarrow^{\tr}&&&&\downarrow^{\tr}&&\downarrow^{\tr}\\
\underline{\ff}_n^{op}&=&\underline{\xx}_{0,n}^{op}&\LA{\Omega}&\underline{\xx}_{1,n-1}^{op}&\LA{\Omega}&\cdots&\LA{\Omega}&\underline{\xx}_{n-1,1}^{op}&\LA{\Omega}&\underline{\xx}_{n,0}^{op}&=&\underline{\ww}_n^{op}
\end{diag}\]

In particular, we have $\ff_n=\Omega ^n \ww_n\subseteq
\Omega^n(\mod\Lambda)$ for any $n\ge 1$.

(2) By (1), $\Omega^n(\mod\Lambda)\subseteq\ff_m$ if and only if
$\tr\Omega^n(\mod\Lambda)\subseteq\ww_m^{op}$ if and only if
$\rgrade{\trf{n+1}{C}}\ge m+1$ holds for any $C\in\mod\Lambda$.

(3) $\ext^1_\Lambda(X,Y)=0$ holds for any $X\in\ww_n$ and $Y\in\pp_{n}$.
In fact, we take a projective resolution
$0\to P_{n-1}\to\cdots\to P_0\to Y\to0$ and apply the functor $\hom_\Lambda(X,\ )$, then we get $\ext^1_\Lambda(X,Y)\simeq\ext^2_\Lambda(X,\Omega Y)\simeq\cdots\simeq \ext^n_\Lambda(X,\Omega^{n-1}Y)=\ext^n_\Lambda(X,P_{n-1})=0$.

(4) Let $C$ be in $\mod\Lambda$. We let $\sigma _{C}: C \rightarrow
C^{**}$ denote the canonical evaluation homomorphism; thus
$\sigma _{C}(x)(f)=f(x)$ for any $x \in C$ and $f \in C^*$. Then there
exists an exact sequence $0\to \extf{1}{\rm Tr}C \to C \buildrel
{\sigma _C} \over \longrightarrow C^{**} \to \extf{2}{\rm Tr}C \to
0$ [AB]. Recall that $C$ is called {\it torsionless} (resp. {\it reflexive}) if
$\sigma _{C}$ is a monomorphism (resp. an isomorphism). So $C$ is
torsionless (resp. reflexive) if and only if it is 1-torsionfree
(resp. 2-torsionfree). In particular, we have $\ff_1=\Omega(\mod\Lambda)$.

(5) For any $X\in\mod\Lambda$ and $n>0$, we have an exact sequence
$0\to\extf{n}{X}\to\trf{n}{X}\stackrel{f}{\to}\Omega\trf{n+1}{X}\to0$
such that $f^*$ is an isomorphism [Ho3].

(6) If $X\in\mod\Lambda$ satisfies $\grade{\extf{i}{X}}>i$ for any
$0\le i\le n$, then $\grade{X}>n$ [Ho2; 6.2][Iy1; 2.3].

(7) For $l,n\ge0$, the following conditions are equivalent [Iy1;
6.1][Hu1; 2.8].

\strut\kern1em(i) $\fd I_i(\Lambda^{op})<l$ holds for any $0 \le
i<n$.

\strut\kern1em(ii) $\sgrade{\extf{l}{C}}\ge n$ holds for any
$C\in\mod\Lambda$.

In this case, we say that $\Lambda^{op}$ satisfies
the {\it $(l,n)$-condition} (or $\Lambda$ satisfies the
{\it $(l,n)^{op}$-condition}). 


\vskip.5em{\bf\XC\ The conditions $G_n(k)$ and $g_n(k)$ }

In this section we introduce Auslander-type conditions $G_n(k)$ and $g_n(k)$ and study their properties. Let us start with the following observation.

\vskip.5em{\bf\XCA\ Lemma }{\it Let $C\in\ff_m$ with $m\ge0$.

(1) $\Omega C\in\ff_{m+1}$ if and only if $\grade{\extf{1}{C}}\ge m$.

(2) $\ee(C,\ff_m)\subseteq\ff_m$ if and only if $\ee(C,\pp_1)\subseteq\ff_m$ if and only if $\sgrade{\extf{1}{C}}\ge m$.}

\vskip.5em{\sc Proof }(1) We have $\tr C\in\ww_m^{op}$ for
$C\in\ff_m$. Since we have an exact sequence
$0\to\extf{1}{C}\to\tr{C}\stackrel{f}{\to}\Omega\tr\Omega C\to0$
such that $f^*$ is an isomorphism by \XBD(5),
$\extf{i}{\extf{1}{C}}=\extf{i+2}{\tr\Omega C}$ holds for any $0 \le
i<m$. Thus the assertion follows.

(2) See [AR3; 1.1].
\hfill\rule{5pt}{10pt}

\vskip.5em{\bf\XCAA\ }
Immediately we have the following result.

\vskip.5em{\bf Lemma }{\it Assume
$\Omega^n(\mod\Lambda)\subseteq\ff_{m}$ with $n,m\ge0$. Then
for the following (1)$\Leftrightarrow$(2)$\Leftarrow$
(3)$\Leftrightarrow$(4) holds.

(1) $\grade{\extf{n+1}{C}}\ge m$ holds for any $C\in\mod\Lambda$.

(2) $\Omega^{n+1}(\mod\Lambda)\subseteq\ff_{m+1}$.

(3) $\sgrade{\extf{n+1}{C}}\ge m$ holds for any $C\in\mod\Lambda$.

(4) $\ee(\Omega^n(\mod\Lambda),\Omega^n(\mod\Lambda))\subseteq\ff_m$.}

\vskip.5em{\sc Proof }
We only have to check (3)$\Leftrightarrow$(4). This is shown by applying \XCA\ as follows:

(3)$\Rightarrow\ee(\Omega^n(\mod\Lambda),\ff_m)\subseteq\ff_m\Rightarrow$(4)$\Rightarrow\ee(\Omega^n(\mod\Lambda),\pp_1)\subseteq\ff_m\Rightarrow$(3).
\hfill\rule{5pt}{10pt}

\vskip.5em{\bf\XCB\ Proposition }{\it
Let $0=l_0<l_1<l_2<l_3<\cdots$ be a (finite or infinite) sequence of integers.
Then for the following (1)$\Leftrightarrow$(2)$\Leftarrow$(3)$\Leftrightarrow$(4)$\Leftrightarrow$(5) holds.

(1) $\grade{\extf{l_i+1}{C}}\ge i$ holds for any $C\in\mod\Lambda$ and $i$.

(2) $\Omega^{l_i+1}(\mod\Lambda)\subseteq\ff_{i+1}$ holds for any $i$.

(3) $\sgrade{\extf{l_i+1}{C}}\ge i$ holds for any $C\in\mod\Lambda$ and $i$.

(4) $\ee(\Omega^{l_i}(\mod\Lambda),\Omega^{l_i}(\mod\Lambda))\subseteq\ff_i$ holds for any $i$.

(5) $\fd I_{i-1}(\Lambda^{op})\le l_i$ holds for any $i$.}

\vskip.5em{\sc Proof } 
(1)$\Leftrightarrow$(2)$\Leftarrow$(3)$\Leftrightarrow$(4)  Since
$\Omega^{l_i}(\mod\Lambda)\subseteq\Omega^{l_{i-1}+1}(\mod\Lambda)$
holds by $l_i\ge l_{i-1}+1$, we can inductively show the assertions by \XCAA.

(3)$\Leftrightarrow$(5) By \XBD(7).
\hfill\rule{5pt}{10pt}

\vskip.5em{\bf\XCBA\ Question }
Put $l_i:=\inf\{ l\ |\ \Omega^{l+1}(\mod\Lambda)\subseteq\ff_{i+1}\}$.
How does the sequence $(l_i)_i$ behave? Is there an example satisfying $l_i=l_{i+1}$?
This equality case was excluded in \XCB.

\vskip.5em{\bf\XCC\ Definition }Let $n,k\ge0$. We say that
$\Lambda$ is $G_n(k)$ if
$\sgrade{\extf{i+k}{C}}\ge i$ holds for any $C\in\mod\Lambda$ and
$1\le i\le n$. By \XBD(7), $\Lambda$ is $G_n(k)$ if and only if $\fd
I_i(\Lambda^{op})\le i+k$ holds for any $0\le i<n$.
Thus $G_n(0)$ is the $n$-Gorenstein
condition, and $G_n(1)^{op}$ is the quasi $n$-Gorenstein condition
in [Hu2] (see \XDB, \XDC\ below).

Similarly, we say that $\Lambda$ is $g_n(k)$
if $\grade{\extf{i+k}{C}}\ge i$ holds for any
$C\in\mod\Lambda$ and $1\le i\le n$. We say that $\Lambda$ is
$G_n(k)^{op}$ (resp. $g_n(k)^{op}$) if $\Lambda^{op}$ is $G_n(k)$
(resp. $g_n(k)$).

We have the following obvious relations for any $n\ge n'$ and $k\le k'$:
\[\begin{diag}
G_{n}(k)&\Rightarrow& G_{n'}(k')\\
\Downarrow&&\Downarrow\\
g_{n}(k)&\Rightarrow& g_{n'}(k')
\end{diag}\]

\vskip0em{\bf\XCD\ Theorem }{\it The conditions (1) and (2) below
are equivalent for any $n,k\ge0$. If $k>0$, then (1)--(5)
are equivalent.

(1) $\Lambda$ is $g_n(k)$.

(2) For any monomorphism $A\stackrel{f}{\to} B$ with
$A,B\in\Omega^{k+1}(\mod\Lambda)$, $\extf{i}{\extf{i}{f}}$ is a
monomorphism for any $0\le i<n$.

(3) $\Omega^{i+k}(\mod\Lambda)\subseteq\ff_{i+1}$ holds for any
$1\le i\le n$.

(4) For any $C\in\mod\Lambda$ and $0 \le i\le n$, there exists an
exact sequence $0\to Y\to X\stackrel{}{\to}\Omega^{k-1}C\to0$ with
$X\in\ww_{i+1}$ and $Y\in\pp_{i+1}$.

(5) For any $C\in\mod\Lambda$ and $0 \le i\le n$, there exists an
exact sequence $0\to\Omega^kC\to Y'\to X'\to0$ with $X'\in\ww_{i+1}$
and $Y'\in\pp_{i+1}$.}

\vskip.5em{\bf\XCDA\ Remark } By \XBD(3), the sequence in (4) gives a
right $\ww_{i+1}$-approximation of $\Omega^{k-1}C$, and the sequence
in (5) gives a left $\pp_{i+1}$-approximation of $\Omega^kC$.

\vskip.5em{\bf\XCDB\ Lemma }{\it Assume that $n \geq 1$ and
$C\in\ww_{n-1}$ satisfies $\grade{\extf{n}{C}}\ge n-1$. Then there
exists an exact sequence $0\to Y\to X\stackrel{}{\to} C\to0$ with
$Y\in\pp_{n}$ and $X\in\ww_n$.}

\vskip.5em{\sc Proof } Let $0\to\Omega^nC\to P_{n-1}\to\cdots\to
P_0\to C\to0$ be a projective resolution of $C$. Take the following
commutative diagram, where the lower sequence is a projective
resolution of $\extf{n}{C}$:
\[\begin{diag}
0&\stackrel{}{\longleftarrow}&\extf{n}{C}&\stackrel{}{\longleftarrow}&
(\Omega^nC)^*&\stackrel{}{\longleftarrow}&P^*_{n-1}&\stackrel{}{\longleftarrow}&
\cdots&\stackrel{}{\longleftarrow}&P^*_1&\stackrel{}{\longleftarrow}&P^*_0\\
&&\parallel&&\uparrow&&\uparrow&&&&\uparrow\\
0&\stackrel{}{\longleftarrow}&\extf{n}{C}&\stackrel{}{\longleftarrow}&Q_0&
\stackrel{}{\longleftarrow}&Q_1&\stackrel{}{\longleftarrow}&
\cdots&\stackrel{}{\longleftarrow}&Q_{n-1}
\end{diag}\]
Taking the mapping cone of the above commutative diagram, we get an
exact sequence $0 \leftarrow (\Omega^nC)^* \leftarrow
P_{n-1}^*\oplus Q_0 \leftarrow \cdots \leftarrow P_0^*\oplus
Q_{n-1}$. Since $\grade{\extf{n}{C}}\ge n-1$, we have an exact sequence
$0\to Q_0^*\to Q_1^*\to\cdots\to Q_{n-1}^*$.
So the last column in the following commutative diagram of
exact sequences is the desired exact sequence:
\[\begin{diag}
&&&&0&&&&0&&0&&0\\
&&&&\downarrow&&&&\downarrow&&\downarrow&&\downarrow\\
&&0&\stackrel{}{\longrightarrow}&Q^*_0&\stackrel{}{\longrightarrow}&\cdots&
\stackrel{}{\longrightarrow}&Q^*_{n-2}&\stackrel{}{\longrightarrow}&Q^*_{n-1}&
\stackrel{}{\longrightarrow}&Y&\stackrel{}{\longrightarrow}&0\\
&&\downarrow&&\downarrow&&&&\downarrow&&\downarrow&&\downarrow\\
0&\stackrel{}{\longrightarrow}&\Omega^nC&\stackrel{}{\longrightarrow}&
P_{n-1}\oplus Q^*_0&\stackrel{}{\longrightarrow}&\cdots&\stackrel{}{\longrightarrow}&
P_1\oplus Q^*_{n-2}&\stackrel{}{\longrightarrow}&
P_0\oplus Q^*_{n-1}&\stackrel{}{\longrightarrow}&X&\stackrel{}{\longrightarrow}&0\\
&&\parallel&&\downarrow&&&&\downarrow&&\downarrow&&\downarrow\\
0&\stackrel{}{\longrightarrow}&\Omega^nC&\stackrel{}{\longrightarrow}&
P_{n-1}&\stackrel{}{\longrightarrow}&\cdots&\stackrel{}{\longrightarrow}&P_1&
\stackrel{}{\longrightarrow}&P_0&\stackrel{}{\longrightarrow}&C&
\stackrel{}{\longrightarrow}&0&\\
&&&&\downarrow&&&&\downarrow&&\downarrow&&\downarrow\\
&&&&0&&&&0&&0&&0
\end{diag}\]
\vskip-2em\hfill\rule{5pt}{10pt}

\vskip0em{\bf\XCDC\ Proof of \XCD\ } (1)$\Rightarrow$(2) For any $A
\in \Omega ^{k+1}(\mod\Lambda)$, we have $A=\Omega ^{k+1}A'$
with $A' \in\mod\Lambda$. So
$\extf{i}{\extf{i}{A}}=\extf{i}{\extf{i+k+1}{A^\prime}}=0$ holds for
any $0<i<n$ by (1). Thus $\extf{i}{\extf{i}{f}}$ is monic. Let us
consider $f^{**}$. If $k>0$, then
$\Omega^{k+1}(\mod\Lambda)\subseteq\ff_2$ holds by \XCAA(1)$\Rightarrow$(2). Thus $A$
and $B$ are reflexive, and $f^{**}\simeq f$ is monic. Let $k=0$. Take
an injection $a:B\to P$ with $P$ projective in $\mod\Lambda$,
and consider an exact sequence $0\to A\stackrel{fa}{\to}P\to
C\to0$ with $C=\Cok(fa)$. Then we have an exact sequence
$P^*\stackrel{(fa)^*}{\longrightarrow}A^*\to\extf{1}{C}\to0$.
Then $f^{**}a^{**}=(fa)^{**}$ is monic by
$\grade{\extf{1}{C}}\ge1$. Thus $f^{**}$ is also monic.

(2)$\Rightarrow$(1) For any $C\in\mod\Lambda$, there exists an exact
sequence $0\to\Omega^{k+1}C\stackrel{f}{\to}P_k\to\Omega^kC\to0$ in
$\mod\Lambda$ with $P_k$ projective. Then we have an exact sequence
$P_k^*\stackrel{f^*}{\to}(\Omega^{k+1}C)^*\to\extf{k+1}{C}\to0$.
Since $f^{**}$ is monic by (2), we obtain $\grade{\extf{k+1}{C}}\ge1$ by
taking $(\ )^*$. On the other hand, since $\extf{i}{\extf{i}{\Omega^{k+1}C}}
\stackrel{\extf{i}{\extf{i}{f}}}{\longrightarrow}\extf{i}{\extf{i}{P_k}}=0$
is monic for any $0<i<n$,
we obtain $0=\extf{i}{\extf{i}{\Omega^{k+1}C}}=\extf{i}{\extf{i+k+1}{C}}$.
Thus $\grade{\extf{i+k+1}{C}}\ge i+1$ holds for
any $C\in\mod\Lambda$ and $0 \le i<n$.

(1)$\Leftrightarrow$(3) Put $l_i:=i+k-1$ in \XCB.

(5)$\Rightarrow$(3) By (5), for any $C\in\mod\Lambda$ and $0 \le
i\le n$, there exists an exact sequence $0\to\Omega^kC\to Y'\to
X'\to0$ in $\mod\Lambda$ with $X'\in\ww_{i+1}$ and $Y'\in\pp_{i+1}$.
Taking $i$-th syzygies, we have an exact sequence $0\to\Omega^{i+k}C\to\Omega^iY'\to\Omega^iX'\to0$ with projective $\Omega^iY'$. Thus we have
$\Omega^{i+k}C=\Omega^{i+1}X'\in\Omega^{i+1}\ww_{i+1}=\ff_{i+1}$ by \XBD(1).

(4)$\Rightarrow$(5) For any $C\in\mod\Lambda$, by (4) there exists
an exact sequence $0\to Y\to X\stackrel{}{\to}\Omega^{k-1}C\to0$ in
$\mod\Lambda$ with $X\in\ww_{i+1}$ and $Y\in\pp_{i+1}$. On the other
hand, there exists an exact sequence $0\to \Omega^kC\to
P_{k-1}\stackrel{}{\to}\Omega^{k-1}C\to0$ with $P_{k-1}$ projective.
Applying \XBA, we have the desired sequence.

(1)$\Rightarrow$(4) We proceed by induction on $i$. 
Assume that $i \geq 0$ and we have an
exact sequence $0\to Y_{i-1}\to X_{i-1}\to\Omega^{k-1}C\to0$ with
$Y_{i-1}\in\pp_{i}$ and $X_{i-1}\in\ww_{i}$. Since
$\extf{i+1}{X_{i-1}}\simeq\extf{i+1}{\Omega^{k-1}C}\simeq\extf{i+k}{C}$
holds, we obtain $\grade{\extf{i+1}{X_{i-1}}}\ge i$ by (1). Applying
\XCDB\ to $X_{i-1}$, we have an exact sequence $0\to Y^\prime\to
X\to X_{i-1}\to0$ with $Y^\prime\in\pp_{i+1}$ and $X\in\ww_{i+1}$.
Taking the following pull-back diagram, the middle row is the
desired exact sequence:

\[\begin{diag}
&&0&&0\\
&&\uparrow&&\uparrow\\
0&\RA{}&Y_{i-1}&\RA{}&X_{i-1}&\RA{}&\Omega ^{k-1}C&\RA{}&0\\
&&\uparrow&&\uparrow&&\parallel\\
0&\RA{}&Y&\RA{}&X&\RA{}&\Omega ^{k-1}C&\RA{}&0\\
&&\uparrow&&\uparrow\\
&&Y^\prime&=&Y^\prime&\\
&&\uparrow&&\uparrow\\
&&0&&0
\end{diag}\]
\vskip-2em\hfill\rule{5pt}{10pt}

\vskip.5em{\bf\XCE\ Theorem }{\it The conditions (1) and (2) below
are equivalent for any $n,k\ge0$. If $k>0$, then (1)--(3)
are equivalent.

(1) $\Lambda$ is $G_n(k)$.

(2) For any exact sequence $0\to A\stackrel{f}{\to} B\to C\to0$ with
$C\in\Omega^k(\mod\Lambda)$, $\extf{i}{\extf{i}{f}}$ is a
monomorphism for any $0 \le i<n$.

(3)
$\ee(\Omega^{i+k}(\mod\Lambda),\Omega^{i+k}(\mod\Lambda))\subseteq\ff_{i+1}$
holds for any $0 \le i<n$.}

\vskip.5em{\sc Proof } (1)$\Rightarrow$(2) Let $C=\Omega^k C^\prime$
with $C' \in {\rm mod}\ \Lambda$ and $g:=\extf{i}{f}$. We have an
exact sequence
$\extf{i+k}{C^\prime}\to\extf{i}{B}\stackrel{g}{\to}\extf{i}{A}\to\extf{i+k+1}{C^\prime}$.
Thus we have exact sequences $\extf{i}{\Cok
g}\to\extf{i}{\extf{i}{A}}\stackrel{a}{\to}\extf{i}{\Im g}$ and
$\extf{i-1}{\Ker g}\to\extf{i}{\Im
g}\stackrel{b}{\to}\extf{i}{\extf{i}{B}}$. Since $\extf{i}{\Cok
g}=0=\extf{i-1}{\Ker g}$ holds by (1), $\extf{i}{g}=ab$ is monic.

(2)$\Rightarrow$(1) For any $C\in {\rm mod}\ \Lambda$, fix $i$ ($0
\le i <n$) and a $\Lambda^{op}$-submodule $D$ of $\extf{i+k+1}{C}$.
Take an exact sequence $Q\stackrel{a}{\to}D\to0$ in mod $\Lambda
^{op}$ with $Q$ projective and $a'$ the composition
$Q\stackrel{a}{\to}D \hookrightarrow \extf{i+k+1}{C}$. We lift $a'$
to $b:Q\to(\Omega^{i+k+1}C)^*$. Take the following push-out diagram,
where $b^\prime$ is the composition $\Omega^{i+k+1}C \buildrel
{\sigma _{\Omega^{i+k+1}C}}\over \longrightarrow
(\Omega^{i+k+1}C)^{**}\buildrel {b^*}\over \longrightarrow Q^*$ and
$P_{i+k}$ is a projective module in $\mod\Lambda$:
\[\begin{diag}
0&\RA{}&\Omega^{i+k+1}C&\RA{}&P_{i+k}&\RA{}&\Omega^{i+k}C&\RA{}&0\\
&&\downarrow^{b^\prime}&&\downarrow&&\parallel\\
0&\RA{}&Q^*&\RA{c}&X&\RA{d}&\Omega^{i+k}C&\RA{}&0
\end{diag}\]
We then have the following commutative diagram with exact rows:
\[\begin{diag}
0&\LA{}&\extf{i+k+1}{C}&\LA{}&(\Omega^{i+k+1}C)^*&\LA{}&
P^*_{i+k}&\LA{}&(\Omega^{i+k}C)^*&\LA{}&0\\
&&\cup&&\uparrow^b&&\uparrow&&\parallel\\
0&\LA{}&D&\LA{a}&Q&\LA{c^*}& X^*&\LA{d^*}&(\Omega^{i+k}C)^*&\LA{}&0
\end{diag}\]

Let $i=0$. Since $c^{**}$ is monic by (2), we obtain $D^*=0$. Thus
$\sgrade{\extf{1+k}{C}}\ge1$.

Fix $i$ ($0<i<n$) and assume that $\Lambda$ is $G_i(k)$. By \XCB,
$\ee(\Omega^{i+k}(\mod\Lambda),\Omega^{i+k}(\mod\Lambda))\subseteq\ff_i$
holds. In particular, $X$ in the first diagram is contained in
$\ff_i$. Take the following commutative diagram with exact rows,
where the upper sequence is still exact by taking $(\ )^*$:
\[\begin{diag}
0&\RA{}&X&\RA{}&Q_0&\RA{}&\cdots&\RA{}&Q_{i-1}&\RA{}&Y&\RA{}&0\\
&&\downarrow^d&&\downarrow&&&&\downarrow&&\downarrow^{g}\\
0&\RA{}&\Omega^{i+k}C&\RA{}&P_{i+k-1}&\RA{}&\cdots&\RA{}&
P_k&\RA{}&\Omega^kC&\RA{}&0\end{diag}\]

We can assume that $g$ is epic by adding a projective direct summand
to $Y$ and $Q_{i-1}$. Thus we obtain an exact sequence $0\to
Z\stackrel{h}{\to}Y\stackrel{g}{\to}\Omega^kC\to0$ such that
$\Omega^ig=d$. Since
$\Cok\extf{i}{h}\simeq\Ker\extf{i+1}{g}\simeq\Ker\extf{1}{d}\simeq\Cok
c^*\simeq D$ and $\extf{i}{Y}=0$ hold, we obtain $D\simeq\extf{i}{Z}$.
Since $\extf{i}{\extf{i}{h}}$ is monic by (2), we have
$\extf{i}{D}=0$. Thus $\sgrade{\extf{i+k+1}{C}}\ge i+1$ holds,
and $\Lambda$ is $G_{i+1}(k)$. Inductively, we obtain (1).

(1)$\Leftrightarrow$(3) Put $l_i:=i+k-1$ in \XCB.\hfill\rule{5pt}{10pt}

\vskip.5em{\bf\XD\ The conditions $G_n(k)$ and $g_n(k)$ for $k=0,1$ }

In this section we concentrate on the conditions $G_n(k)$ and $g_n(k)$ for the case $k=0,1$. Let us start with giving a quick proof of the following
remarkable `left-right symmetry',
where (1) is well-known in [FGR; 3.7], (2) is in [HN;
4.7][Hu3; 2,4] and (3) is in [AR3][HN; 4.1].

\vskip.5em{\bf\XDA\ Theorem }{\it
(1) $G_n(0)\Leftrightarrow G_n(0)^{op}$.\ \ (2) $g_n(1)\Leftrightarrow g_n(1)^{op}$.\ \
(3) $g_n(0)\Leftrightarrow G_n(1)^{op}$.}

\vskip.5em{\sc Proof }We shall proceed by using induction on $n$.
The case $n=0$ is obvious. Now assume that $n \ge 1$ and the assertions hold
for $n-1$. Thus we can assume that $g_{n-1}(1)^{op}$ holds for each
case, so $\Omega^i(\mod\Lambda^{op})=\ff_i^{op}$ holds for any $1
\le i\le n$ by \XBD(1)(4) and \XCD. Thus $\trf{i+1}{D}\in\ww_{i}$ holds for any
$D\in\mod\Lambda^{op}$ and $1\le i\le n$.

(1) We will show the `only if' part. Take $D\subseteq\extf{n}{C}$ for
$C\in\mod\Lambda^{op}$. We have an exact sequence $0\to
D\to\trf{n}{C}\stackrel{f}{\to}D^\prime\to0$ such that $f^*$ is an
isomorphism by \XBD(5). Since $\trf{n}{C}\in\ww_{n-1}$ holds, we
have $\extf{i}{D}\subseteq\extf{i+1}{D^\prime}$ for any $0\le i<n$.
Since $\grade{\extf{i}{D}}>i$ holds for any $0\le i<n$ by $G_n(0)$,
we have $\grade{D}\ge n$ by \XBD(6).

(2) We will show the `only if' part. Put $D:=\extf{n+1}{C}$ for $C\in\mod\Lambda^{op}$.
We have an exact sequence $0\to
D\to\trf{n+1}{C}\stackrel{f}{\to}\Omega\trf{n+2}{C}\to0$ such that
$f^*$ is an isomorphism by \XBD(5). Since $\trf{n+1}{C}\in\ww_{n}$
holds, we have
$\extf{i}{D}\simeq\extf{i+1}{\Omega\trf{n+2}{C}}\simeq\extf{i+2}{\trf{n+2}{C}}$
for any $0\le i<n$. Since $\grade{\extf{i}{D}}>i$ holds for any
$0\le i<n$ by $g_n(1)$, we have $\grade{D}\ge n$ by \XBD(6).

(3) We will show the `only if' part. Take $D\subseteq\extf{n+1}{C}$ for
$C\in\mod\Lambda^{op}$. We have an exact sequence $0\to
D\to\trf{n+1}{C}\stackrel{f}{\to}D^\prime\to0$ such that $f^*$ is an
isomorphism by \XBD(5). Since $\trf{n+1}{C}\in\ww_{n}$ holds, we
have $\extf{i}{D}\simeq\extf{i+1}{D^\prime}$ for any $0\le i<n$.
Since $\grade{\extf{i}{D}}>i$ holds for any $0\le i<n$ by $g_n(0)$,
we have $\grade{D}\ge n$ holds by \XBD(6).

We will show the `if' part. Put $D:=\extf{n}{C}$ for
$C\in\mod\Lambda$. We have an exact sequence $0\to
D\to\trf{n}{C}\stackrel{f}{\to}\Omega\trf{n+1}{C}\to0$ such that
$f^*$ is an isomorphism by \XBD(5). Since $\trf{n}{C}\in\ww_{n-1}^{op}$
holds, we have
$\extf{i}{D}\subseteq\extf{i+1}{\Omega\trf{n+1}{C}}\simeq
\extf{i+2}{\trf{n+1}{C}}$ for any $0\le i<n$. Since
$\grade{\extf{i}{D}}>i$ holds for any $0\le i<n$ by $G_n(1)^{op}$, we
have $\grade{D}\ge n$ by \XBD(6).\hfill\rule{5pt}{10pt}

\vskip.5em{\bf\XDAA\ Question } 
It is natural to ask for the existence of a common generalization
of the conditions $G_n(k)$ and $g_n(k)$ satisfying 
certain `left-right symmetry'. For example,
is there some natural condition $G_n(k,l)$ for each triple $(n,k,l)$
of non-negative integers with the following properties?

(i) $G_n(k,0)=G_n(k)$, and $G_n(k,1)=g_n(k)$.

(ii) $G_n(k,l)\Leftrightarrow G_n(l,k)^{op}$.

(iii) $G_n(k,l)\Rightarrow G_{n^\prime}(k^\prime,l^\prime)$
if $n\ge n^\prime$, $k\le k^\prime$ and $l\le l^\prime$.

\vskip.5em{\bf\XDB\ Theorem }(cf. [FGR; 3.7]) {\it The following
conditions are equivalent.

(1) $\Lambda$ is $G_n(0)$, i.e. $n$-Gorenstein.

(2) $\extf{i}{\extf{i}{}}$ preserves monomorphisms in $\mod\Lambda$
for any $0\le i<n$.

(3) $\fd I_i(\Lambda)\le i$ holds for any $0\le i<n$.

($i$)$^{op}$ Opposite side version of ($i$) ($1\leq i \leq 3$).}

\vskip.5em{\bf\XDC\ Theorem }(cf. [AB][Ho1; 2.1][HN; 4.7][Hu2;
3.3][Hu3; 2.4]) {\it The following conditions are equivalent.

(1) $\Lambda$ is $g_n(1)$.

(2) For any monomorphism $A\stackrel{f}{\to} B$ with
$A,B\in\Omega^{2}(\mod\Lambda)$, $\extf{i}{\extf{i}{f}}$ is a
monomorphism for any $0\le i<n$.

(3) $\Omega^i(\mod\Lambda)=\ff_i$ holds for any $1\le i \leq n+1$.

(4) For any $C\in\mod\Lambda$ and $0\le i \leq n$, there exists an
exact sequence $0\to Y\to X\to C\to0$ with $X\in\ww_{i+1}$ and
$Y\in\pp_{i+1}$.

(5) For any $C\in\mod\Lambda$ and $0\le i \leq n$, there exists an
exact sequence $0\to\Omega C\to Y\to X\to0$ with $X\in\ww_{i+1}$ and
$Y\in\pp_{i+1}$.

($i$)$^{op}$ Opposite side version of ($i$) ($1\le i\le 5$).}

\vskip.5em{\bf\XDD\ Theorem }(cf. [Ho1; 2.4][AR3][IST; 2.1][HN;
4.1][Hu2; 3.6]) {\it The following conditions (1)--(7) are
equivalent. If $\Lambda$ is an artin algebra, then (8) is also
equivalent.

(1) $\Lambda$ is $G_n(1)$.

(2) For any exact sequence $0\to A\stackrel{f}{\to} B\to C\to0$ with
$C\in\Omega(\mod\Lambda)$, $\extf{i}{\extf{i}{f}}$ is a monomorphism for any
$0 \le i<n$.

(3)
$\ee(\Omega^{i+1}(\mod\Lambda),\Omega^{i+1}(\mod\Lambda))\subseteq\ff_{i+1}$
holds for any $0 \le i<n$.

(4) $\Omega^i(\mod\Lambda)$ is closed under extensions for any $1 \le i\le n$.

(4$'$) $\add\Omega^i(\mod\Lambda)$ is closed under extensions for any
$1\le i\le n$.

(5) $\fd I_{i}(\Lambda^{op})\le i+1$ holds for any $0 \le i<n$.

(6) $\Lambda$ is $g_n(0)^{op}$.

(7) For any monomorphism $A\stackrel{f}{\to} B$ with
$A,B\in\Omega(\mod\Lambda^{op})$, $\extf{i}{\extf{i}{f}}$ is a
monomorphism for any $0\le i<n$.

(8) For any $C\in\mod\Lambda$ and $1\le i\le n$, there exist exact
sequences $0\to Y\to X\stackrel{}{\to} C\to0$ and $0\to C\to
Y^\prime\to X^\prime\to0$ with $X,X^\prime\in\Omega^i(\mod\Lambda)$
and $Y,Y^\prime\in\ii_{i+1}$.}

\vskip.5em{\bf\XDE\ Question } (1) Is it possible to characterize
$n$-Gorenstein rings in terms of the categories
$\Omega^i(\mod\Lambda)$ (like \XDC(3), \XDD(4))
or the existence of approximation sequences (like \XDC(4)(5), \XDD(8))?

(2) When does the equivalence with \XDD(8) hold for noetherian
rings? When is the category $\Omega^i(\mod\Lambda)$
contravariantly finite for a noetherian ring $\Lambda$? In this case,
we have a sequence in \XDD(8) by Wakamatsu's Lemma \XBCA.

\vskip.5em{\bf\XDF\ Proof of \XDB--\XDD\ } \XDB, \XDC\ and
(1)$\Leftrightarrow$(2)$\Leftrightarrow$(3)$\Leftrightarrow$(5)$\Leftrightarrow$(6)$\Leftrightarrow$(7)
in \XDD\ follow immediately from \XCB, \XCD, \XCE\ and \XDA.
We now show the other implications in \XDD.

(3)$\Rightarrow$(4)$\Rightarrow$(4$'$) Easy.

(4$'$)$\Rightarrow$(3) Assume $\Omega^i(\mod\Lambda)=\ff_i$ for some $i\le n$.
Then $\add\Omega^i(\mod\Lambda)=\ff_i$ holds, and we have
$\ee(\Omega^i(\mod\Lambda),\Omega^i(\mod\Lambda))\subseteq\ff_i$ by (4$'$).
By \XCAA(4)$\Rightarrow$(2), we have $\Omega^{i+1}(\mod\Lambda)=\ff_{i+1}$.
Thus we have (3) inductively.

(8)$\Rightarrow$(4$'$) Let $0\to A\stackrel{}{\to}B\to C\to 0$ be
an exact sequence with $A,C\in\add\Omega^i(\mod\Lambda)$. Take
an exact sequence $0\to Y\to X\stackrel{f}{\to}B\to 0$ with
$X\in\Omega^i(\mod\Lambda)$ and $Y\in\ii_{i+1}$. Since we have an exact sequence
$0=\ext^1_\Lambda(C,Y)\to\ext^1_\Lambda(B,Y)\to\ext^1_\Lambda(A,Y)=0$,
we have that $f$ splits. Thus $B\in\add\Omega^i(\mod\Lambda)$ holds.


(4)+(5)$\Rightarrow$(8) For any $C\in\mod\Lambda$, take the
following commutative diagram whose upper sequence is the minimal
injective resolution of $C$ and vertical sequences are minimal
projective resolutions: {\small\[\begin{diag}
&&&&&&0&&0&&&&0\\
&&&&&&\uparrow&&\uparrow&&&&\uparrow\\
0&\longrightarrow&C&\longrightarrow&I_0&\longrightarrow&I_1&\longrightarrow&I_2&
\longrightarrow&\cdots&\longrightarrow&I_i\\
&&&&&&\uparrow&&\uparrow&&&&\uparrow\\
&&&&&&P_{1,0}&\longrightarrow&P_{2,0}&\longrightarrow&\cdots&\longrightarrow&P_{i,0}\\
&&&&&&&&\uparrow&&&&\uparrow\\
&&&&&&&&P_{2,1}&\longrightarrow&\cdots&\longrightarrow&P_{i,1}\\
&&&&&&&&&&&&\uparrow\\
&&&&&&&&&&&&\vdots\\
&&&&&&&&&&&&\uparrow\\
&&&&&&&&&&&&P_{i,i-1}
\end{diag}\]}

Then we obtain the following commutative diagram of exact sequences:

{\small\[\begin{diag}
&&0&&0&&&&0&&0&&0\\
&&\uparrow&&\uparrow&&&&\uparrow&&\uparrow&&\uparrow\\
0&\longrightarrow&X&\longrightarrow&\bigoplus_{j=0}^{i-1}P_{j+1,j}&
\longrightarrow&\cdots&\longrightarrow&P_{i-1,0}\oplus
P_{i,1}&\longrightarrow&
P_{i,0}&\longrightarrow&\Omega^{-i-1}C&\longrightarrow&0\\
&&\uparrow&&\uparrow&&&&\uparrow&&\uparrow&&\uparrow\\
0&\longrightarrow&Y&\longrightarrow&I_0\oplus\bigoplus_{j=0}^{i-1}P_{j+1,j}&
\longrightarrow&\cdots&\longrightarrow&I_{i-2}\oplus P_{i-1,0}\oplus P_{i,1}&
\longrightarrow&I_{i-1}\oplus P_{i,0}&\longrightarrow&I_i&\longrightarrow&0\\
&&\uparrow&&\uparrow&&&&\uparrow&&\uparrow&&\uparrow\\
0&\longrightarrow&C&\longrightarrow&I_0&\longrightarrow&\cdots&\longrightarrow&I_{i-2}&
\longrightarrow&I_{i-1}&\longrightarrow&\Omega^{-i}C&\longrightarrow&0\\
&&\uparrow&&\uparrow&&&&\uparrow&&\uparrow&&\uparrow\\
&&0&&0&&&&0&&0&&0
\end{diag}\]} where $\Omega^{-i-1}C={\rm
Cok}(I_{i-1}\to I_i)$ and $\Omega^{-i}C={\rm Cok}(I_{i-2}\to I_{i-1})$. Since
$\id(I_l\oplus\bigoplus_{j=0}^{i-l-1}P_{j+l+1,j})\le i-l$ for any $0
\le l \le i-1$ holds by (5), we have $Y\in\ii_{i+1}$. Thus we get one
of the desired sequences $0\to C\to Y\to X\to 0$.

Now we take exact sequences $0\to\Omega C\to P\to C\to0$ with projective $P$
and $0\to\Omega C\to Y^\prime\to X^\prime\to0$ with
$X^\prime\in\Omega^i(\mod\Lambda)$ and $Y^\prime\in\ii_{i+1}$.
Since $\Omega^i(\mod\Lambda)$ is closed under extensions by (4),
we have the desired sequence $0\to Y^\prime\to
X^{\prime\prime}\to C\to0$ by applying the dual of \XBA.
\hfill\rule{5pt}{10pt}

\vskip.5em{\bf\XDG\ Corollary }{\it (1) If $\Lambda$ is $g_n(1)$, then $\pp_{i+1}$ is a covariantly finite
subcategory of $\mod\Lambda$ for any $0 \le i\le n$.

(2)[Hu2; 3.6] If an artin algebra $\Lambda$ is $g_n(0)$, then
$\pp_{i+1}$ is a functorially finite subcategory of $\mod\Lambda$ for
any $0 \le i\le n$.}

\vskip.5em{\sc Proof } (1) For any $C\in\mod\Lambda$, let $C^\prime$
be a maximal factor module of $C$ such that any homomorphism from
$C$ to $\pp_{i+1}$ factors through $C^\prime$. We can take an exact sequence
$0\to C^\prime\to Y\to C^{\prime\prime}\to0$ with $Y\in\pp_{i+1}$. Take
the sequence $0\to Y^\prime\to X^\prime\to C^{\prime\prime}\to0$
with $X^\prime\in\ww_{i+1}$ and $Y^\prime\in\pp_{i+1}$ by \XDC(4). 
Applying \XBA, we have an exact sequence
$0\to C'\stackrel{f}{\to}Y''\to X'\to0$ with $Y''\in\pp_{i+1}$.
Then $f$ is a left $\pp_{i+1}$-approximation of $C'$, and
so the composition $C\to
C^{\prime}\stackrel{f}{\to}Y^{\prime\prime}$ gives a left
$\pp_{i+1}$-approximation of $C$.

(2) $\pp_{i+1}$ is covariantly finite by (1). Since $\ii_{i+1}^{op}$ is
covariantly finite by \XDD(8), $\pp_{i+1}$ is contravariantly finite.
\hfill\rule{5pt}{10pt}

\vskip.5em{\bf\XDH\ Theorem }{\it Assume that an artin algebra
$\Lambda$ is $G_n(1)$. For any $i\ge0$ and $n+1\ge j\ge1$ satisfying
$i+j\le n+2$,
$(\xx_{i,j-1},\yy_{i,j})$ forms a functorially finite cotorsion pair.}

\vskip.5em{\sc Proof }Since $\xx_{i,j-1}$
and $\yy_{i,j}$ are closed under direct summands, we only have to show
that $(\xx_{i,j-1},\yy_{i,j})$ satisfies the conditions \XBC(3)(4) by
\XBCC. $\ext^1_\Lambda(\ww_{i},\pp_{i})=0$ by \XBD(3), and
$\ext^1_\Lambda(\ff_{j-1},\ii_{j})=0$ by $\ff_{j-1}\subseteq
\Omega^{j-1}(\mod\Lambda)$,
so $\ext^1_\Lambda(\xx_{i,j-1},\yy_{i,j})=0$ holds.
For any $C\in\mod\Lambda$, we can take an exact sequence
$0\to I\to \Omega^{j-1}C^\prime\to C\to0$ with $I\in\ii_{j}$
and $C'\in\mod\Lambda$ by \XDD(8). By \XDC(4),
we can take an exact sequence $0\to P\to X\to C^\prime\to0$ with
$X\in\ww_{i+j-1}$ and $P\in\pp_{i+j-1}$. Taking $(j-1)$-th syzygies, we
have an exact sequence $0\to P^\prime\to
X^\prime\to\Omega^{j-1}C^\prime\to0$ with $X^\prime\in\xx_{i,j-1}$ and
$P^\prime\in\pp_{i}$. Consider the following pull-back diagram:

\[\begin{diag}
&&0&&0&&\\
&&\uparrow&&\uparrow&&\\
0&\RA{}&I&\RA{}&\Omega^{j-1}C^\prime&\RA{}&C&\RA{}&0\\
&&\uparrow&&\uparrow&&\parallel\\
0&\RA{}&Y&\RA{}&X^\prime&\RA{}&C&\RA{}&0\\
&&\uparrow&&\uparrow\\
&&P^\prime&=&P^\prime\\
&&\uparrow&&\uparrow&&\\
&&0&&0&&
\end{diag}\]
Since the middle row $0\to Y\to X^\prime\to C\to0$ satisfies
$Y\in\yy_{i,j}$ and $X^\prime\in\xx_{i,j-1}$, the condition \XBC(3) holds.

Moreover, take an exact sequence $0\to C\to I\to\Omega^-C\to0$ with injective
$I$ and $0\to Y_0\to X_0\to\Omega^-C\to0$ with $X_0\in\xx_{i,j-1}$ and
$Y_0\in\yy_{i,j}$. Applying the proof of \XBA\ and using
$\ee(\ii_1,\yy_{i,j})\subseteq\yy_{i,j}$, we have an exact sequence
$0\to C\to Y_1\to X_1\to0$ with $X_1\in\xx_{i,j-1}$ and $Y_1\in\yy_{i,j}$.
Thus the condition \XBC(4) holds.\hfill\rule{5pt}{10pt}

\vskip.5em{\bf \XDI\ }Our main theorem \XAC\ is a special case of the following result.

\vskip.5em{\bf Corollary }{\it Let $\Lambda$ be an artin algebra
which is $G_\infty(1)$ and $G_\infty(1)^{op}$. Then
$(\xx_{i,j-1},\yy_{i,j})$ ($i\ge0,\ j\ge1$) and
$(\yy_{i,j},\ddd\xx_{j,i-1}^{op})$ ($i\ge1,\ j\ge0$) form
functorially finite cotorsion pairs.}

\vskip.5em{\sc Proof } This is immediate from 4.8 and the fact that
$\yy_{i,j}=\ddd\yy_{j,i}^{op}$.\hfill\rule{5pt}{10pt}

\vskip.5em{\bf\XDJ\ }Denote by $\ww_{\infty}:=\bigcap _{n \geq 1}\ww_n
=\{ C\in\mod\Lambda\ |\ \ext^i_\Lambda(C,\Lambda)=0 \mbox{ for any } 
i \geq 1\}$. As an application of \XAC, we have the following

\vskip.5em{\bf Corollary }{\it Let $\Lambda$ be an
$\infty$-Gorenstein artin algebra. Then $\ww_1 \supseteq \ww_2
\supseteq \cdots \supseteq \ww_{\infty}$ is a chain of
contravariantly finite subcategories of $\mod\Lambda$.}

\vskip.5em{\sc Proof }By \XAC, we have that $\ww_i$ is
contravariantly finite for any $i \geq 1$. By [AR1; 6.12] and [AR2;
5.5(b)], $\ww_{\infty}$ is contravariantly finite.
\hfill\rule{5pt}{10pt}

\vskip.5em{\bf\XE\ Finitistic dimension and the $(l,n)$-condition }

In this final section, we give some results on finitistic dimension and
left-right symmetry of the $(l,n)$-condition.

\vskip.5em{\bf\XEA\ }Recall that the finitistic dimension of
$\Lambda$, denoted by $\fin\Lambda$, is defined as $\sup\{\pd
X\ |\ X\in\mod\Lambda \mbox{ and } \pd X<\infty\}$.

\vskip.5em{\bf Lemma }{\it Assume that $\Lambda$ is $g_{n+1}(k)$ with
$n \geq 0$ and $k>0$. If $\fin\Lambda =n$, then id $\Lambda \leq n+k$.}

\vskip.5em{\sc Proof }
Let $C\in\mod\Lambda$. By \XCD, there exists an exact
sequence $0 \to Y \to X \to \Omega ^{k-1}C \to 0$ with $X\in \ww
_{n+2}$ and $\pd Y \leq n+1$. Then $\pd Y \leq n$ because
fin.dim$\Lambda =n$. So $\extf{n+k+1}{C}\simeq \extf{n+2}{\Omega
^{k-1}C}\simeq \extf{n+1}{Y}=0$ and hence $\id \Lambda \leq
n+k$.\hfill\rule{5pt}{10pt}

\vskip.5em{\bf\XEB\ Theorem }{\it If $\Lambda$ is $g_\infty(k)$ with
$k\geq 0$, then $\fin\Lambda \leq \id \Lambda\leq\fin\Lambda +k$.}

\vskip.5em{\sc Proof }
It is well known that $\fin\Lambda \leq \id \Lambda$.
So it suffices to prove $\id\Lambda \leq\fin\Lambda +k$.
The case $k>0$ follows easily from \XEA.
Now suppose $k=0$. Then $\Lambda$ {\it is} $g_\infty(0)$, and so
$\fd I_i(\Lambda) \leq i+1$ for any $i \geq 0$ by [HN; 4.1] (see \XDD\ above).
It follows from [Hu3; 2.15] that $\fin\Lambda=\id \Lambda$.
\hfill\rule{5pt}{10pt}

\vskip.5em{\bf\XEC\ Corollary }{\it
(1) If $\Lambda$ is $G_\infty(0)$, then $\fin\Lambda=\id \Lambda$
and $\fin\Lambda^{op}=\id \Lambda ^{op}$.

(2) If $\Lambda$ is $G_\infty(1)$, then
$\fin\Lambda \leq \id \Lambda \leq\fin\Lambda+1$ and
$\fin\Lambda ^{op}=\id \Lambda ^{op}$.}

\vskip.5em{\sc Proof }
(1) follows from the symmetry of $G_\infty(0)$ and \XEB.

(2) Because $\Lambda$ is $G_\infty(1)$ if and only if $\Lambda$ is
$g_\infty(0)^{op}$ by [AR3; 0.1] and [HN; 4.1] (see \XDD\ above),
our conclusion follows from \XEB\ and its dual.\hfill\rule{5pt}{10pt}

\vskip.5em{\bf\XED\ }Recall that we say that $\Lambda$ satisfies
the {\it $(l,n)$-condition} (or $\Lambda^{op}$ satisfies
the {\it $(l,n)^{op}$-condition}) if $\sgrade{\extf{l}{C}}\ge n$ holds
for any $C\in\mod\Lambda^{op}$ (see \XBD(7)). Similarly,
we say that $\Lambda$ satisfies
the {\it weak $(l,n)$-condition} (or $\Lambda^{op}$ satisfies
the {\it weak $(l,n)^{op}$-condition}) if $\grade{\extf{l}{C}}\ge n$ holds
for any $C\in\mod\Lambda^{op}$.
For example, $\Lambda$ is $G_n(k)$ if and only if $\Lambda$ satisfies the
$(k+i,i)^{op}$-condition for any $1 \leq i \leq n$, and $\Lambda$ is
$g_n(k)$ if and only if $\Lambda$ satisfies the weak
$(k+i,i)^{op}$-condition for any $1 \leq i \leq n$.

\vskip.5em{\bf Lemma }{\it
(1) $(k,l)\ +$ weak $(l,n)\Rightarrow$ $(k,n)$.

(2) weak $(k,l)\ +$ weak $(l,n)\Rightarrow$ weak $(k,n)$.

(3) $(k,l)\ +$ weak $(l,n)^{op}\Rightarrow$ $(k,n)$.

(4) weak $(k,l)\ +$ weak $(l,n)^{op}\Rightarrow$ weak $(k,n)$.}

\vskip.5em{\sc Proof }Using [Iy1; 2.3(2)], we can show all assertions
in a manner similar to the proof of [Iy1; 2.3(3)].
\hfill\rule{5pt}{10pt}

\vskip.5em{\bf\XEE\ }We call $l\in\nnn$ a {\it dominant number} of
$\Lambda$ if $\fd I_i(\Lambda)<\fd I_l(\Lambda)$ holds for any $i<l$
[Iy1]. The following result generalizes a theorem in [Iy1; 1.1].

\vskip.5em{\bf Theorem }{\it Any dominant number $l$ of
$\Lambda$ satisfies $\fd I_l(\Lambda)\ge l$.}

\vskip.5em{\sc Proof }Put $m:=\fd I_l(\Lambda)$ and assume $m<l$.
Since $l$ is a dominant number of $\Lambda$, $\Lambda$ satisfies the
$(m,l)$ and $(m+1,l+1)$-condition. Since $\Lambda$ satisfies the
$(m,m+1)$-condition by $m<l$, $\Lambda$ satisfies the
$(m,l+1)$-condition by the dual of \XED(1). This is in contradiction to
$m=\fd I_l(\Lambda)$.\hfill\rule{5pt}{10pt}

\vskip.5em{\bf\XEF\ }The following result is an analog of \XED(2) and \XDA.

\vskip.5em{\bf Theorem }{\it Let $k,l,n\ge0$. Assume $k\ge n-1$.

(1) If $\sgrade{\extf{1}{\ff_k}}\ge l$ and
$\sgrade{\extf{1}{\ff_l^{op}}}\ge n$, then
$\sgrade{\extf{1}{\ff_k}}\ge n$.

(2) If $\grade{\extf{1}{\ff_k}}\ge l-1$ and
$\grade{\extf{1}{\ff_l^{op}}}\ge n-1$, then
$\grade{\extf{1}{\ff_k}}\ge n-1$.}

\vskip.5em{\sc Proof } We can assume $n>l$.

(1) Fix $C\in\ff_k$ and $D\subseteq\extf{1}{C}$. Then $\grade{D}\ge
l$ holds by our first assumption. By \XBD(5), we have an exact sequence $0\to
D\to\tr C\stackrel{f}{\to}D^\prime\to0$ such that $\tr C\in\ww_k^{op}$
and $f^*$ is an isomorphism.
Thus $D^\prime\in\ww_l^{op}$ by $k\ge l$, and $\Omega^lD^\prime\in\ff_l^{op}$.
Using \XCA(1) and $\sgrade{\extf{1}{\ff_l^{op}}}\ge n>l$,
we have $\Omega^iD^\prime\in\ff_l^{op}$ for any $i\ge l$ inductively.
On the other hand, applying $(\ )^*$ to the exact sequence above,
we have an inclusion $\extf{i}{D}\subseteq
\extf{i+1}{D^\prime}=\extf{1}{\Omega^iD^\prime}$ for any $0\le i\le k$.
Thus our second assumption implies $\grade{\extf{i}{D}}\ge n$ for any
$l\le i\le k$. Since $\grade{D}\ge l$, we obtain $\grade{D}\ge n$ by \XBD(6).

(2) Fix $C\in\ff_k$ and put $D:=\extf{1}{C}$ and
$D^\prime:=\Omega\tr\Omega C$. Then $\grade{D}\ge l-1$ and $\Omega C\in\ff_l$
hold by \XCA(1). We have $\tr\Omega C\in\ww_l^{op}$ and
$\Omega^{l-1}D^\prime=\Omega^l\tr\Omega C\in\ff_l^{op}$.
Using \XCA(1) and $\grade{\extf{1}{\ff_l^{op}}}\ge n-1>l-1$,
we have $\Omega^iD^\prime\in\ff_l^{op}$ for any $i\ge l-1$ inductively.
By \XBD(5), we have an exact sequence $0\to D\to\tr C\stackrel{f}{\to}
D^\prime\to0$ such that $\tr C\in\ww_k^{op}$ and $f^*$ is an isomorphism.
Applying $(\ )^*$, we have
$\extf{i}{D}\simeq\extf{i+1}{D^\prime}\simeq\extf{1}{\Omega^iD^\prime}$
for any $0\le i<k$.
Thus our second assumption implies $\grade{\extf{i}{D}}\ge n-1$ for
any $l-1\le i<k$. Since $\grade{D}\ge l-1$, we obtain
$\grade{D}\ge n-1$ by \XBD(6).\hfill\rule{5pt}{10pt}

\vskip.5em{\bf\XEG\ } Recall that $\Lambda$ is $G_n(0)$ if and only if
$\Lambda$ satisfies the $(i,i)^{op}$-condition for any $1 \leq i \leq n$.
It is known that $\Lambda$ is $G_n(0)$ if and only if so is $\Lambda
^{op}$ [FGR; 3.7] (see \XDB\ above).
However, the $(i, i)$-condition doesn't possess
such a symmetric property in general [Hu1; p.1460][Iy2; 2.1.1].
For example, the algebra given by the quiver $\bullet\stackrel{a}{\to}\bullet
\stackrel{b}{\to}\bullet\leftarrow\bullet$ modulo the ideal $ab$
satisfies exactly one of the $(2,2)$ and $(2,2)^{op}$-conditions.

As an application of \XEF, we give a sufficient condition that the
$(i,i)$-condition implies the $(i,i)^{op}$-condition as follows.

\vskip.5em{\bf Corollary } {\it $G_{n-1}(1)+(n,n)\Rightarrow(n,n)^{op}$.}

\vskip.5em{\sc Proof }The case $n=0$ is trivial.
Now suppose $n \geq 1$. Then
$\sgrade{\extf{1}{\ff_{n-1}}}\ge n-1$ and
$\sgrade{\extf{1}{\ff_{n-1}^{op}}}\ge n$. So
$\sgrade{\extf{1}{\ff_{n-1}}}\ge n$ by \XEF(1). On the other hand,
$\Omega^i(\mod\Lambda)=\ff_i$ for any $1 \le i \le n$ by [HN; 4.7]
(see \XCD\ above). Thus $\sgrade{\extf{n}C}=\sgrade{\extf{1}{\Omega
^{n-1}C}}\geq n$ for any $C\in\mod\Lambda$ and $\Lambda$ satisfies
the $(n,n)^{op}$-condition.\hfill\rule{5pt}{10pt}

\vskip.5em{\bf\XEGA\ Corollary }{\it $(2,2)^{op}+(3,3)\Rightarrow(3,3)^{op}$.}

\vskip.5em{\sc Proof }This is immediate from \XEG\ for
$n=3$.\hfill\rule{5pt}{10pt}

\vskip.5em{\bf\XEGB\ Example } 
We give an example satisfying the
conditions in \XEGA. Let $K$ be a field and $\Lambda$ a finite
dimensional $K$-algebra given by the quiver
\[\def\arraystretch{.5}\begin{array}{c}
\bullet\\
\downarrow\\
\bullet\stackrel{a}{\longrightarrow}\bullet\longrightarrow\bullet\\
\ \downarrow^b\\
\bullet\end{array}\]
modulo the ideal $ab$. Then $\fd I_0(\Lambda)=\fd I_1(\Lambda)=
\fd I_0(\Lambda ^{op})=\fd I_1(\Lambda^{op})=1$, and
$\fd I_2(\Lambda)=\fd I_2(\Lambda ^{op})=2$.

\vskip.5em{\bf\XEH\ Corollary }{\it The following conditions are equivalent.

(1) $\sgrade{\extf{i}{C}}\ge i$ holds for any $C\in\mod\Lambda$ and
$1\le i\le n$ (i.e. $\Lambda$ is $G_n(0)$).

(2) $\sgrade{\extf{1}{\ff_{i-1}}}\ge i$ holds for any $1 \le i\le
n$.

(i)$^{op}$ The opposite side version of ($i$) ($i=1,2$).}

\vskip.5em{\sc Proof } (1)$\Rightarrow$(2) is trivial.

(2)$\Rightarrow$(1) By \XCA(1), $\Omega\ff_{i-1}\subseteq\ff_i$ holds
for any $1\le i\le n$. This implies inductively that
$\Omega^i(\mod\Lambda)=\ff_i$ for any $1 \le i\le n$. So the
assertion follows easily.

(2)$^{op}\Rightarrow$(2) The case $n=0$ is obvious.
Now suppose $n\ge 1$. By \XEF(1),
$\sgrade{\extf{1}{\ff_{n-1}}}\ge n-1$ and
$\sgrade{\extf{1}{\ff_{n-1}^{op}}}\ge n$ imply
$\sgrade{\extf{1}{\ff_{n-1}}}\ge n$. Thus the assertion follows
inductively.\hfill\rule{5pt}{10pt}

\vskip.5em{\small
\vskip.5em{\bf References }

[A] M. Auslander: Representation dimension of Artin algebras. Lecture notes,
Queen Mary College, London, 1971.

[AB] M. Auslander, M. Bridger: Stable module theory. Memoirs of the
American Mathematical Society,
No. 94 American Mathematical Society, Providence, R.I. 1969.

[AR1] M. Auslander, I. Reiten: Applications of contravariantly
finite subcategories. Adv. Math. 86 (1991), no. 1, 111--152. 

[AR2] M. Auslander, I. Reiten: $k$-Gorenstein algebras and syzygy
modules. J. Pure Appl. Algebra 92 (1994), no. 1, 1--27.

[AR3] M. Auslander, I. Reiten: Syzygy modules for Noetherian rings.
J. Algebra 183 (1996), no. 1, 167--185.

[B] J.-E. Bjork: The Auslander condition on Noetherian rings. Seminaire d'Algebre Paul Dubreil et Marie-Paul Malliavin, 39eme Annee (Paris, 1987/1988), 137--173, Lecture Notes in Math., 1404, Springer, Berlin, 1989. 

[C] J. Clark: Auslander-Gorenstein rings for beginners. International Symposium on Ring Theory (Kyongju, 1999), 95--115, Trends Math., Birkhauser Boston, Boston, MA, 2001.

[EJ] E. E. Enochs, O. M. G. Jenda: Gorenstein balance of Hom and tensor.  Tsukuba J. Math.  19  (1995),  no. 1, 1--13. 

[EX] E. E. Enochs, J. Xu: Gorenstein flat covers of modules over Gorenstein rings. J. Algebra 181 (1996), no. 1, 288--313. 

[EHIS] K. Erdmann, T. Holm, O. Iyama, J. Schr\"oer: Radical embeddings and representation dimension, Adv. Math. 185 (2004), no. 1, 159--177.

[FGR] R. M. Fossum, P. Griffith, I. Reiten: Trivial extensions of abelian categories.
Homological algebra of trivial extensions of abelian categories with applications to ring theory.
Lecture Notes in Mathematics, Vol. 456. Springer-Verlag, Berlin-New York, 1975.

[FI] K. R. Fuller, Y. Iwanaga: On $n$-Gorenstein rings and Auslander rings of low injective dimension. Representations of algebras (Ottawa, ON, 1992), 175--183, CMS Conf. Proc., 14, Amer. Math. Soc., Providence, RI, 1993.

[Ho1] M. Hoshino: Syzygies and Gorenstein rings. Arch. Math. (Basel) 55 (1990), no. 4, 355--360.


[Ho2] M. Hoshino: Noetherian rings of self-injective dimension two,
Comm. Algebra, 21 (1993), no. 4, 1071--1094.

[Ho3] M. Hoshino: On self-injective dimension of artinian rings,
Tsukuba J. Math, 18 (1994), no. 1, 1--8.


[HN] M. Hoshino, K. Nishida: A generalization of the Auslander
formula, in: Representations of Algebras and Related Topics, Fields
Institute Communications 45, Amer. Math. Soc., Providence, Rhode
Island, 2005, pp.175--186.

[Hu1] Z. Huang: Extension closure of $k$-torsionfree modules, Comm.
Algebra, 27 (1999), no. 3, 1457-1464.

[Hu2] Z. Huang: Approximation presentations of modules and
homological conjectures, to appear in Comm. Algebra.

[Hu3] Z. Huang: On the grade of modules over noetherian rings, to
appear in Comm. Algebra.

[IST] K. Igusa, S. O. Smal$\phi$, G. Todorov: Finite projectivity
and contravariant finiteness, Proc. Amer. Math. Soc. 109 (1990), no.
4, 937--941.

[IS] Y. Iwanaga, H. Sato: On Auslander's $n$-Gorenstein rings. J. Pure Appl. Algebra 106 (1996), no. 1, 61--76.

[Iw1] Y. Iwanaga: On rings with self-injective dimension$\leq 1$.  Osaka J. Math.  15  (1978), no. 1, 33--46.

[Iw2] Y. Iwanaga: On rings with finite self-injective dimension. II.  Tsukuba J. Math.  4  (1980), no. 1, 107--113.

[Iy1] O. Iyama: Symmetry and duality on $n$-Gorenstein rings,
J. Algebra 269 (2003), no. 2, 528--535.

[Iy2] O. Iyama: The relationship between homological properties and representation theoretic realization of Artin algebras. Trans. Amer. Math. Soc. 357 (2005), no. 2, 709--734.

[Iy3] O. Iyama: $\tau$-categories III. Auslander orders and
Auslander-Reiten quivers. Algebr. Represent. Theory 8 (2005), no. 5, 601--619.

[Iy4] O. Iyama: Finiteness of Representation dimension, Proc. Amer. Math. Soc. 131 (2003), no. 4, 1011--1014.

[Iy5] O. Iyama: Higher-dimensional Auslander-Reiten theory on maximal orthogonal subcategories.  Adv. Math.  210  (2007),  no. 1, 22--50.

[Iy6] O. Iyama: Auslander correspondence, Adv. Math. 210 (2007), no. 1, 51--82.

[M] J. Miyachi: Injective resolutions of Noetherian rings and cogenerators. Proc. Amer. Math. Soc. 128 (2000), no. 8, 2233--2242.

[R] R. Rouquier: Representation dimension of exterior algebras. Invent. Math. 165 (2006), no. 2, 357--367.

[Sa] L. Salce: Cotorsion theories for abelian groups. Symposia Mathematica, Vol. XXIII (Conf. Abelian Groups and their Relationship to the Theory of Modules, INDAM, Rome, 1977), pp. 11--32, Academic Press, London-New York, 1979. 

[Sm] S. P. Smith: Some finite-dimensional algebras related to elliptic curves. Representation theory of algebras and related topics (Mexico City, 1994), 315--348, CMS Conf. Proc., 19, Amer. Math. Soc., Providence, RI, 1996.

[T] H. Tachikawa: Quasi-Frobenius rings and generalizations.
Lecture Notes in Mathematics, Vol. 351. Springer-Verlag, Berlin-New York, 1973.

[W] T. Wakamatsu: Tilting modules and Auslander's Gorenstein property. J. Algebra 275 (2004), no. 1, 3--39.
}
\end{document}